\documentclass[preprint, a4paper, 3p]{elsarticle}

\usepackage{graphicx}
\usepackage{amsmath}
\usepackage{amsfonts}
\usepackage{amssymb}
\usepackage{subfig}
\usepackage{xcolor}
\usepackage{multirow}

\newcommand{\vekt}[1]{\mathbf{#1}}
\newcommand{\diff}[1]{\nu \vekt{\nabla}^2 \vekt{v} }
 
\newcommand{\vel}[0]{\vekt{u}^f}
\newcommand{\xtil}[0]{\vekt{\tilde{x}}}
\newcommand{\CauchyStress}[0]{\vekt{T}}
\newcommand{\iterS}[1]{\color{brown}{$#1$}}
\newcommand{\iterF}[1]{\color{blue}{$#1$}}
\newcommand{\iterC}[1]{$\mathbf{#1}$}
\newcommand{\iterN}[1]{$\underline{#1}$}
\newcommand{\iterCompare}[1]{$\textcolor{gray}{#1}$}

\journal{IX International Conference on Computational Methods for Coupled Problems in Science and Engineering}
\journal{conference proceedings of COUPLED PROBLEMS 2021}

\begin{document}
\begin{frontmatter}

\title{THE PERFORMANCE IMPACT OF NEWTON ITERATIONS PER SOLVER CALL IN PARTITIONED FLUID-STRUCTURE INTERACTION}

\author{THOMAS SPENKE$^{*}$, NORBERT HOSTERS$^{*}$ AND MAREK BEHR$^{*}$}

\address{$^{*}$Chair for Computational Analysis of Technical Systems (CATS)\\
	Center for Simulation and Data Science (JARA-CSD) \\
RWTH Aachen University\\
Schinkelstraße 2, 52062 Aachen, Germany\\
e-mail: \{spenke,hosters,behr\}@cats.rwth-aachen.de, \\
web page: http://www.cats.rwth-aachen.de}
%
%


\begin{abstract}
	The cost of a partitioned fluid-structure interaction scheme is typically assessed by the number of coupling iterations required per time step,
while ignoring the Newton loops within the nonlinear sub-solvers.
In this work, we discuss why these single-field iterations 
deserve more attention when evaluating the coupling's efficiency
and how to find the optimal number of Newton steps per coupling iteration.
\end{abstract}

\end{frontmatter}

\section{INTRODUCTION}

Partitioned algorithms enjoy great popularity in the fluid-structure interaction (FSI) community: Treated as black boxes, the fluid and structural solvers are coupled only via the exchange of interface data. 
Over the last decade, the inherent drawbacks regarding stability were mitigated significantly, e.g., by interface quasi-Newton methods \cite{lindner2015comparison,Degroote2009,spenke2020multi}. 

Assuming the cost for data exchange to be negligible, the efficiency of a partitioned scheme is typically assessed by the number of coupling iterations, i.e., solver calls, required per time step. When coupling two nonlinear solvers, however, the cost of one solver call is not constant, but depends on various parameters. In particular, it is closely connected to the number of Newton iterations performed. 
On the one hand, this supports the conclusion that the Newton iterations required per time step provide a much better measure for the coupling's efficiency. On the other hand, it raises a central question: How many Newton iterations should best be run per coupling step?
The answer is not trivial: While always iterating to full convergence produces unnecessary overhead, running too few steps brings the risk of impeding stability by feeding back defective data into the coupling loop. 

This work discusses the impact of the Newton iterations per solver call 
and proposes an adaptive choice to improve the efficiency of partitioned fluid-structure interaction schemes.

\section{PARTITIONED FLUID-STRUCTURE INTERACTION} \label{Section2}

Although the key aspects of this work are expected to hold true for other coupled multi-physics simulations as well,
we restrict ourselves to partitioned fluid-structure interaction.
More precisely, we consider an incompressible fluid in the domain $\Omega^f$ interacting with an elasto-dynamic structure $\Omega^s$. 

\subsection{Fluid Subproblem} \label{Section2.1}
The fluid velocity $\vekt{u}^f(\vekt{x},t)$ and its pressure $p^f(\vekt{x},t)$ are governed by the unsteady Navier-Stokes equations for an incompressible fluid:
\begin{subequations}
	\begin{alignat}{2}
		\rho^f \left( \frac{\partial \vel}{\partial t} + \vel \cdot \boldsymbol{\nabla} \vel - \vekt{b}^f \right) - \boldsymbol{\nabla}  \cdot \CauchyStress^f &= \vekt{0}	\qquad	&& \text{in} ~\Omega_t^f ~~\forall t \in [0,T] ~,\\
		\boldsymbol{\nabla}  \cdot \vel \,&= 0 && \text{in} ~\Omega_t^f ~~\forall t \in [0,T]~,
	\end{alignat}
\end{subequations}
where $\rho^f$ is the constant fluid density, while $\vekt{b}^f$ denotes the resultant of all external body forces per unit mass of fluid.
Assuming a Newtonian fluid 
with dynamic viscosity $\mu^f$,
Stokes' law models the Cauchy stress tensor as
$\CauchyStress^f( \vel, p^f) = - p^f \vekt{I} + \mu^f \left(  \nabla \vel + (\nabla \vel)^T \right).$
The problem is closed by a divergence-free initial velocity field as well as an appropriate set of boundary conditions on $\Gamma^{f}=\partial \Omega^f$.

The fluid problem is simulated by our in-house solver XNS, using stabilized P1P1 finite elements in space \cite{pauli2016stabilized,donea2003finite} and a BDF1 scheme in time \cite{forti2015semi}.
The ALE mesh is adapted to deforming domains via the linear elastic mesh-update method(EMUM) \cite{emumPaper}.

\subsection{Structural Subproblem}  \label{Section2.2}

The structural displacement $\vekt{d}^s(\vekt{x},t)$ is given by the dynamic balance of stresses. From a Lagrangian viewpoint with respect to the undeformed configuration $\Omega_0^s$, it reads
\begin{alignat}{2}
	\rho^s \frac{d^2 \vekt{d}^s}{dt^2} &= \boldsymbol{\nabla}_0 \cdot \left( \vekt{S} \vekt{F}^T \right) + \vekt{b}^s \qquad &&\text{in } \Omega_0^s ~~\forall t \in [0,T]~,
\end{alignat}
where $\rho^s$ denotes the material density and $\vekt{b}^s$ the resultant of all body forces per unit volume. 
The 2nd Piola-Kirchhoff stresses $\vekt{S}$ are defined based on the Cauchy stress tensor $\CauchyStress^s$ and the deformation gradient $\vekt{F}$ 
as $\vekt{S} := \det (\vekt{F}) ~ \vekt{F}^{-1}  \, \CauchyStress^s  \, \vekt{F}^{-T}$.

As constitutive equation, the St. Venant-Kirchhoff material model provides the linear stress-strain law $\vekt{S} = \vekt{C}^{s}:\vekt{E}$, 
where the constant matrix $ \vekt{C}^{s}$ depends on two material parameters, e.g., Young's modulus $E^s$ and Poisson's ratio $\nu^s$.
The  definition of  the Green-Lagrange strains $\vekt{E} := \frac{1}{2} \left( \vekt{F}^T \vekt{F} - \vekt{I} \right)$ 
introduces a geometrical nonlinearity into the structural model \cite{bathe2006finite, yibin2001nonlinear}.

A closed problem formulation requires both an initial displacement field (typically zero) and a set of Dirichlet and Neumann boundary conditions.
The structural problem is then solved by the in-house finite-element code FEAFA using Lagrangian finite elements or isogeometric analysis (IGA) \cite{cottrell2009isogeometric,Hughes} in space and a generalized-$\alpha$ scheme in time \cite{chung1993time,erlicher2002analysis}.

\subsection{Coupling Conditions} \label{Section2.3}

Naturally, in fluid-structure interaction the solution fields are not independent, but instead connected at the shared interface 
$\Gamma^{fs} = \partial \Omega^f \cap \partial \Omega^s$ \cite{Hosters2017}:
\begin{enumerate}
	\item The \textit{kinematic} coupling condition states the continuity of displacements, i.e.,  $\vekt{d}^f = \vekt{d}^s ~\text{on } \Gamma^{fs}$,
	which directly implies the equality of velocities and accelerations, too.
	\item Following Newton's third law, the \textit{dynamic} condition requires the equality of interface tractions:
	$\CauchyStress^f  ~ \vekt{n}^f  = \CauchyStress^s  ~ \vekt{n}^s  ~\text{on } \Gamma^{fs}$, where  $\vekt{n}^f$ and $\vekt{n}^s$ are the normal vectors.
\end{enumerate}

Satisfying these coupling conditions for every point in time,  i.e., in a continuous manner, ensures the conservation of mass, momentum, and energy
over the FSI boundary \cite{kuttler2006solution}.

\subsection{Dirichlet-Neumann Scheme}

This work relies on a partitioned FSI algorithm, meaning the 
two subproblems are addressed by two distinct solvers, that are coupled only via the exchange of interface
data. While this strategy features a high flexibility 
regarding the solvers,
their communication requires some additional considerations:
(1) Since the meshes in general do not match at the interface, a conservative projection is needed, the \textit{spatial coupling} \cite{Hosters2017}.
(2) The interdependency of the two subproblems requires an iterative procedure 
to find a consistent solution of the coupled problem, referred to as \textit{temporal coupling} \cite{degroote2011multi,hostersspline}.

The most common temporal coupling algorithm for FSI problems is the \textit{Dirichlet-Neumann scheme}:
While the fluid tractions are passed as a Neumann boundary condition to the structure (dynamic continuity),
the resulting interface deformation poses a Dirichlet condition for the fluid velocity (kinematic continuity).
For each time step, the two solvers are successively called in a Gauss-Seidel iteration until convergence is reached \cite{hostersspline}. \\

As a partitioned FSI algorithm, the Dirichlet-Neumann scheme suffers from an inherent instability,
caused by the \textit{added-mass effect} \cite{forster2007robust,forster2007artificial,Causin}.
Basically, it is characterized by overestimated deformations causing exaggerated fluid inertia terms and vice versa.

A common countermeasure is to augment the Dirichlet-Neumann scheme by an update step:
Typically, the computed interface deformation $\xtil^k$ is modified before it is passed back to the fluid as $\vekt{x}^{k+1}$.
The simplest version is a \textit{relaxation}, i.e., 
$\vekt{x}^{k+1} =  \omega \xtil^k + (1-\omega) \vekt{x}^k$,
with the relaxation factor $\omega$ being either some constant $\omega < 1$ (``under-relaxation")
or updated dynamically, e.g., via \textit{Aitken's relaxation} \cite{kuttler2008fixed,irons1969version}.

A more sophisticated approach are \textit{interface quasi-Newton (IQN)} methods, which use a Newton-like update
based on an approximated Jacobian of the coupled problem.
The Jacobian approximation is successively
improved by collecting information from the intermediate results of each coupling iteration.
This way, IQN methods almost completely overcome the added-mass difficulty \cite{lindner2015comparison,Degroote2009,spenke2020multi}.

\section{Newton Iterations per Solver Call}

A central aspect of partitioned schemes 
for fluid-structure interaction 
is the concept of treating the single-field solvers
as black boxes, in that only their in- and output, but not their interior setup and techniques, are known. 
Combined with the idea that the solver calls are the most expensive part of the simulation, 
increasing the algorithm's efficiency comes down to decreasing the number of coupling iterations required for convergence.
This implicitly assumes the cost of one solver call to be constant - which is far from true for nonlinear solvers, 
but the best we can do without knowing any internal details.

\subsection{Newton Iterations and Computational Cost} \label{Sec:NewtonCost}

In practice, however, this full black-box case is rather uncommon,
as the user or developer typically has access to either the solver's input configuration, such as the number of Newton iterations performed per  call, or even to the source code.
This naturally raises questions about the Newton iterations' effect on the computational cost.
%

Of course it is impossible to take all effects into account, but
one Newton iteration basically corresponds to one repetition of the numerical solution procedure, i.e., assembling the matrix-vector system and solving it for the Newton increment.
%
%
This supports the conclusion 
that assuming the cost of a solver call to scale with the number of Newton iterations performed
is more accurate than assuming constant cost.

With that, the cost of the simulation 
depends not only on the total number of coupling steps $N_{Coupling}$, but also
the Newton iterations $N_{Newton}^i$ run for each subproblem $i$:
\begin{align*}
	\mathbf{cost}(\textit{simulation}) \approx N_{Coupling} \cdot \mathbf{cost}(\textit{data transfer}) + \sum_{i=f,s} N_{Newton}^i \cdot \mathbf{cost}(\textit{Newton}^i)~.
\end{align*}
As one Newton iteration is typically much more expensive than the data transfer within one coupling iteration,
i.e., $\mathbf{cost}(\textit{Newton}^i) >> \mathbf{cost}(\textit{data transfer})$,
 this work uses the sum of Newton iterations 
$N_{Newton}=\sum_{i=f,s} N_{Newton}^i$
as the main efficiency measure for the partitioned scheme.
Based on that, 
focus is put on how it is influenced by the number of Newton steps run per solver call.
The relation is non-trivial as running too few iterations brings the risk of feeding inaccurate data back into the coupling loop, while
with too many Newton iterations computational time is misspend on polishing up a solution that will be overwritten in the next coupling step anyway.

The remainder of this chapter is outlined as follows:
Section \ref{Sec:Convergence} defines the convergence criteria used in this work,
before Section \ref{Sec:BottomCase} and \ref{Sec:ChannelCase} discuss the effect of the Newton iterations per solver call for two numerical examples.
Drawing conclusions from this, 
Section \ref{Sec:AdaptiveMethods} proposes new approaches for choosing the number of Newton iterations dynamically
and investigates their effectiveness.

\subsection{Convergence Criteria} \label{Sec:Convergence}

Before looking into numerical examples,
it is essential to define the convergence criteria:
\begin{enumerate}
	\item A subproblem is considered converged when its residual vanishes, i.e., is lower than the bound $\varepsilon_{Problem}$. In that case, the problem's Newton loop terminates independently from the iteration number.
	We will refer to this as \textit{single-field convergence}. 
	\item In line with the fixed-point character of the partition scheme \cite{lindner2015comparison,kuttler2008fixed}, \textit{Coupling convergence}
	is reached if the solutions of all subproblems stay virtually unchanged within one coupling iteration,
	 i.e., all relative changes are lower than $\varepsilon_{Coupling}$.
\end{enumerate}

If and only if both conditions are satisfied, the coupled simulation has converged and proceeds with the next time step.

\subsection{Example: Tank with Elastic Bottom} \label{Sec:BottomCase}

\begin{figure}[h!]
	\centering
	\subfloat[Test case illustration.] 	
	{
		\resizebox{0.44\textwidth}{!}{	\input{figures/BottomCase.pdf_tex}} 	\label{Fig:BottomCaseA}
	} 
\qquad
	\subfloat[Snapshot in deformed state ($t=0.2\,\text{s}$).]
	{
		\includegraphics[trim=0 120 0 200,clip, width=0.44\textwidth]{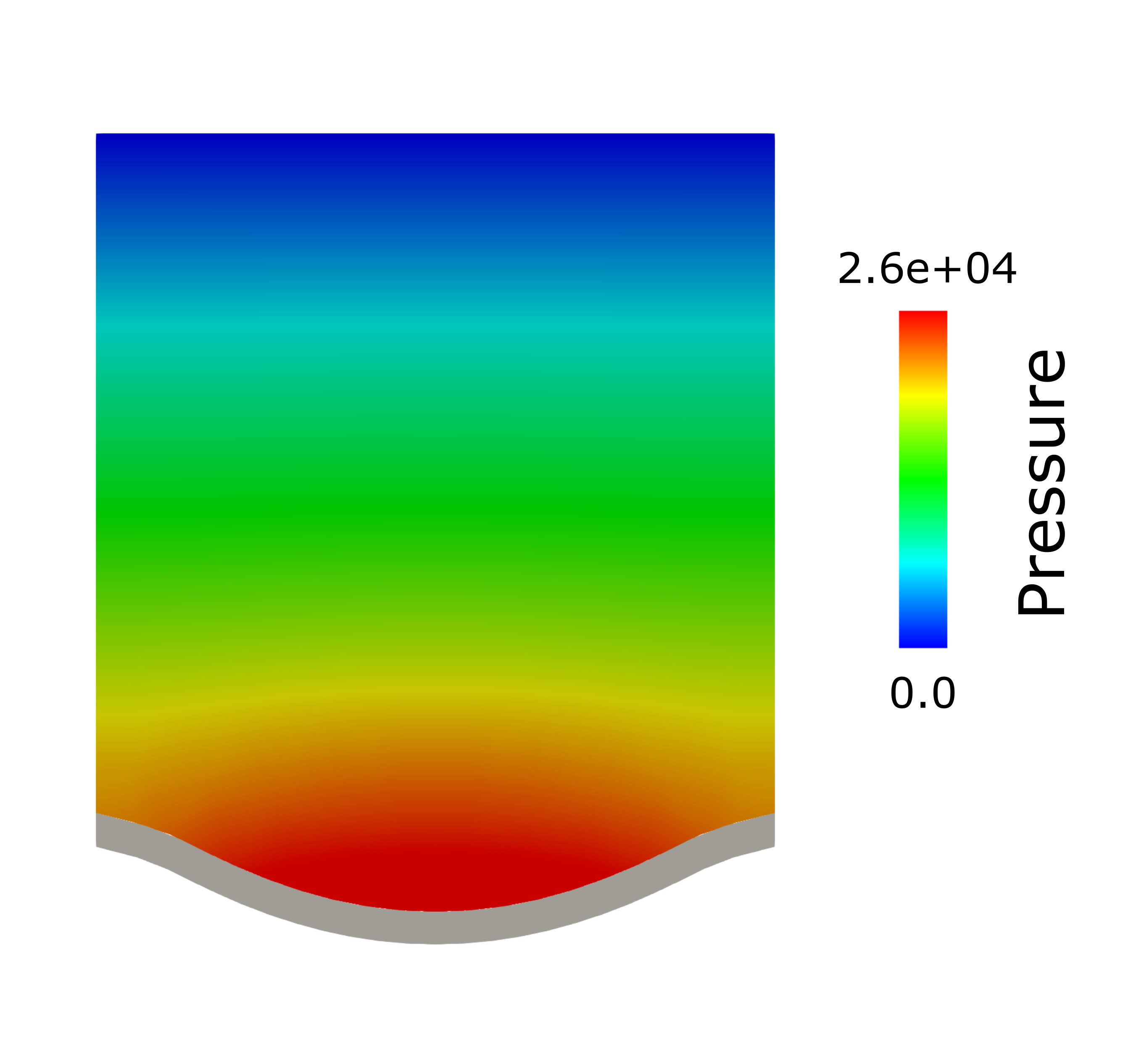}
	}  
	\caption{Elastic bottom test case.}
	\label{Fig:BottomCase}
\end{figure}

The first numerical example of this work is depicted in Figure \ref{Fig:BottomCase}: A square tank with rigid side walls is filled by an incompressible fluid.
For simplicity, the fluid domain has a natural Neumann boundary on the top rather than a free surface, allowing for a free in- and outflow.
Its elastic bottom is clamped on both ends and 
deforms due to the fluid's gravity.
The material and geometrical properties are listed in Figure \ref{Fig:BottomCaseA}.

To simplify the parameter study, the discretization is rather coarse:
While the fluid mesh has $400$ triangular finite elements,
the structural problem is solved by $60$ isogeometric elements of spline degree $2$.
The simulations are run for $50$ time steps of size $\Delta t=0.01$.

To handle the added-mass instability, which is expected to be strong because of the high density ratio $\rho^f / \rho^s = 1.0$, 
the structural deformation is updated via the IQN-IMVLS method \cite{spenke2020multi}.
The convergence criteria were chosen as $\varepsilon_{Coupling}=10^{-5}$ for the coupling and $\varepsilon_{Problem}=10^{-10}$ for both fluid and structure.\\

This test case was run for different combinations of the Newton iterations per solver call for the fluid and the structural subproblem, from now on referred to as $N^f$ and $N^s$.

\begin{table}
		\caption{Iterations required for different choices of $N^f$ and $N^s$. The bold numbers denote the \textbf{coupling iterations}, the underlined ones the \underline{total Newton iterations}. The single-field Newton iterations are colored  \textcolor{blue}{blue for the fluid} and  \textcolor{brown}{brown for the structural solver}.}
	\centering
\small
	\begin{tabular}{ c   c  || c c | c  c | c c | c c | c c |  c c }
		\multicolumn{2}{ c || }{ } &  \multicolumn{12}{ c  } { \textbf{Structural Newton Iterations $N^s$} } \\
		& & \multicolumn{2}{| c |}{\textbf{1}} & \multicolumn{2}{| c |}{\textbf{2}} & \multicolumn{2}{| c |}{\textbf{3}} & \multicolumn{2}{| c |}{\textbf{4}} & \multicolumn{2}{| c |}{\textbf{5}}& \multicolumn{2}{| c }{$\mathbf{\infty}$} \\
		\hline \hline 
		\multirow{1}{*}{ \rotatebox{90}{ \textbf{Fluid Newton Iterations $N^f$} }}
		& \multirow{2}{*}{\textbf{1}}    &  \iterC{617} & \iterN{1166} & \iterC{567} & \iterN{1316}  &  \iterC{583} & \iterN{1499}  &  \iterC{605} & \iterN{1622} &  \iterC{573} & \iterN{1546} &  \iterC{573} & \iterN{1546}  \\
		& 	&   \iterF{$617$} &  \iterS{$567$} & \iterF{$567$} &  \iterS{$549$} & \iterF{$583$} &  \iterS{$549$} & \iterF{$605$} &  \iterS{$549$} & \iterF{$573$} &  \iterS{$549$} & \iterF{$573$} &  \iterS{$549$} \\
				\cline{2-14}
		& \multirow{2}{*}{\textbf{2}}    &  \iterC{497} & \iterN{1488} & \iterC{460} & \iterN{1613}  &  \iterC{438} & \iterN{1667}  &  \iterC{460} & \iterN{1836} &  \iterC{461} & \iterN{1837} &  \iterC{461} & \iterN{1837}  \\
& 	&   \iterF{$994$} &  \iterS{$494$} & \iterF{$919$} &  \iterS{$694$} & \iterF{875} &  \iterS{$792$} & \iterF{919} &  \iterS{$917$} & \iterF{922} &  \iterS{$915$} & \iterF{922} &  \iterS{$915$} \\
		\cline{2-14}
& \multirow{2}{*}{\textbf{3}}    &  \iterC{468} & \iterN{1872} & \iterC{366} & \iterN{1708}  &  \iterC{366} & \iterN{1843}  &  \iterC{374} & \iterN{1974} &  \iterC{381} & \iterN{2014} &  \iterC{381} & \iterN{2016}  \\
& 	&   \iterF{1404} &  \iterS{468} & \iterF{1097} &  \iterS{611} & \iterF{1098} &  \iterS{745} & \iterF{1122} &  \iterS{852} & \iterF{1143} &  \iterS{871} & \iterF{1143} &  \iterS{873} \\
		\cline{2-14}
& \multirow{2}{*}{\textbf{4}}    &  \iterC{407} & \iterN{2011} & \iterC{342} & \iterN{1939}  &  \iterC{335} & \iterN{2043}  &  \iterC{352} & \iterN{2225} &  \iterC{352} & \iterN{2240} &  \iterC{352} & \iterN{2240}  \\
& 	&   \iterF{1604} &  \iterS{407} & \iterF{1342} &  \iterS{597} & \iterF{1315} &  \iterS{728} & \iterF{1386} &  \iterS{839} & \iterF{1386} &  \iterS{854} & \iterF{1386} &  \iterS{854} \\
		\cline{2-14}
& \multirow{2}{*}{\textbf{5}}    &  \iterC{414} & \iterN{2300} & \iterC{335} & \iterN{2108}  &  \iterC{336} & \iterN{2264}  &  \iterC{352} & \iterN{2441} &  \iterC{355} & \iterN{2460} &  \iterC{355} & \iterN{2460}  \\
& 	&   \iterF{1886} &  \iterS{414} & \iterF{1523} &  \iterS{585} & \iterF{1528} &  \iterS{736} & \iterF{1602} &  \iterS{839} & \iterF{1604} &  \iterS{856} & \iterF{1604} &  \iterS{856} \\
	\cline{2-14}
%
%
& \multirow{2}{*}{$\mathbf{\infty}$}    &  \iterC{412} & \iterN{2410} & \iterC{343} & \iterN{2237}  &  \iterC{333} & \iterN{2333}  &  \iterC{355} & \iterN{2565} &  \iterC{354} & \iterN{2547} &  \iterC{354} & \iterN{2547}  \\
& 	&   \iterF{1998} &  \iterS{412} & \iterF{1639} &  \iterS{589} & \iterF{1609} &  \iterS{724} & \iterF{1713} &  \iterS{852} & \iterF{1700} &  \iterS{847} & \iterF{1700} &  \iterS{847} \\
	\end{tabular}
	\label{Tab:BottomCase-Newton} 
\end{table}

\begin{figure}[h!]
	\centering
	\subfloat[Coupling iterations $N_{Coupling}$.]
	{
		\includegraphics[trim=0 10 0 5,clip, width=0.48\textwidth]{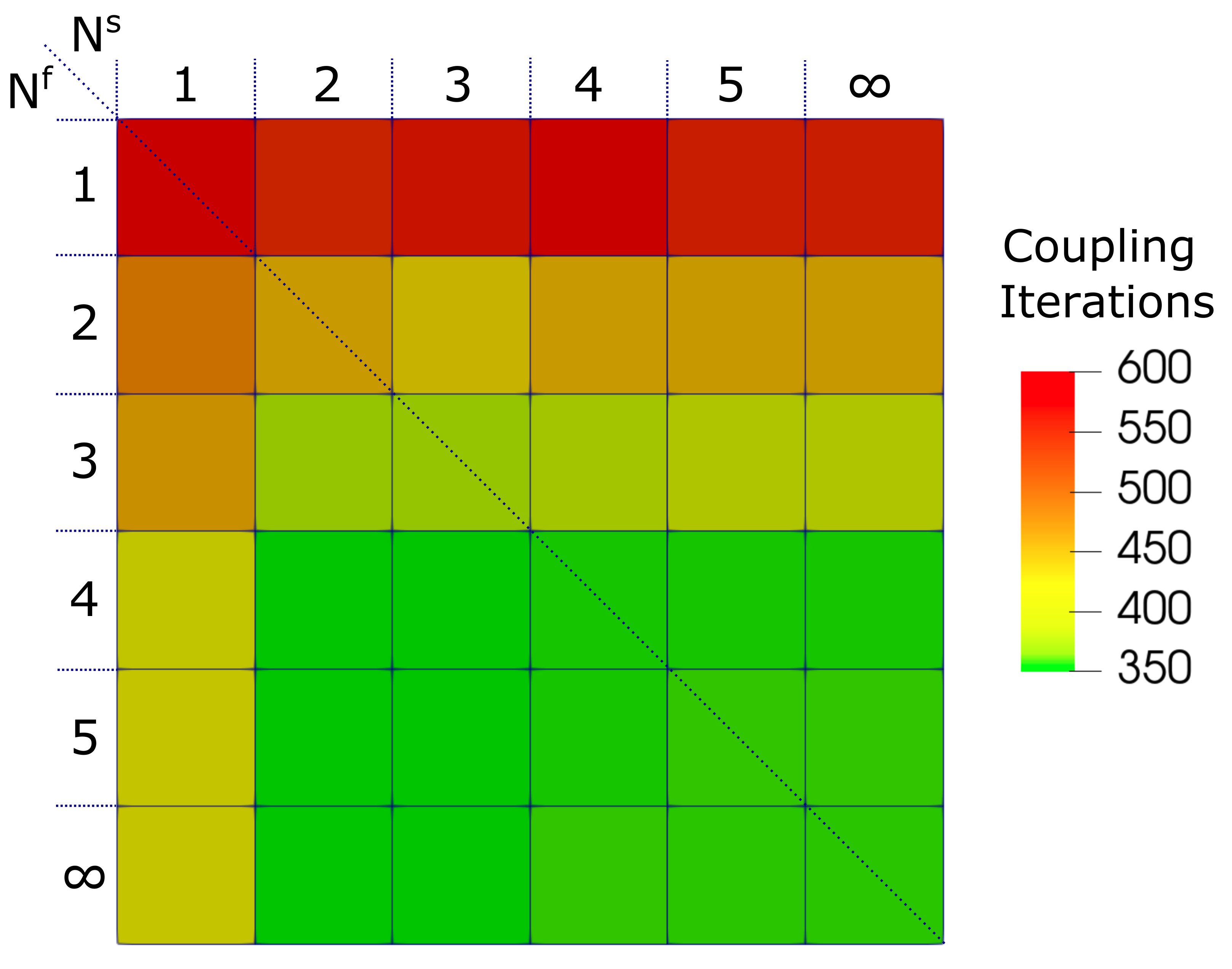}
	}  	
	\subfloat[Total Newton iterations $N_{Newton}$.]
	{
	 	\includegraphics[trim=0 10 0 5,clip, width=0.48\textwidth]{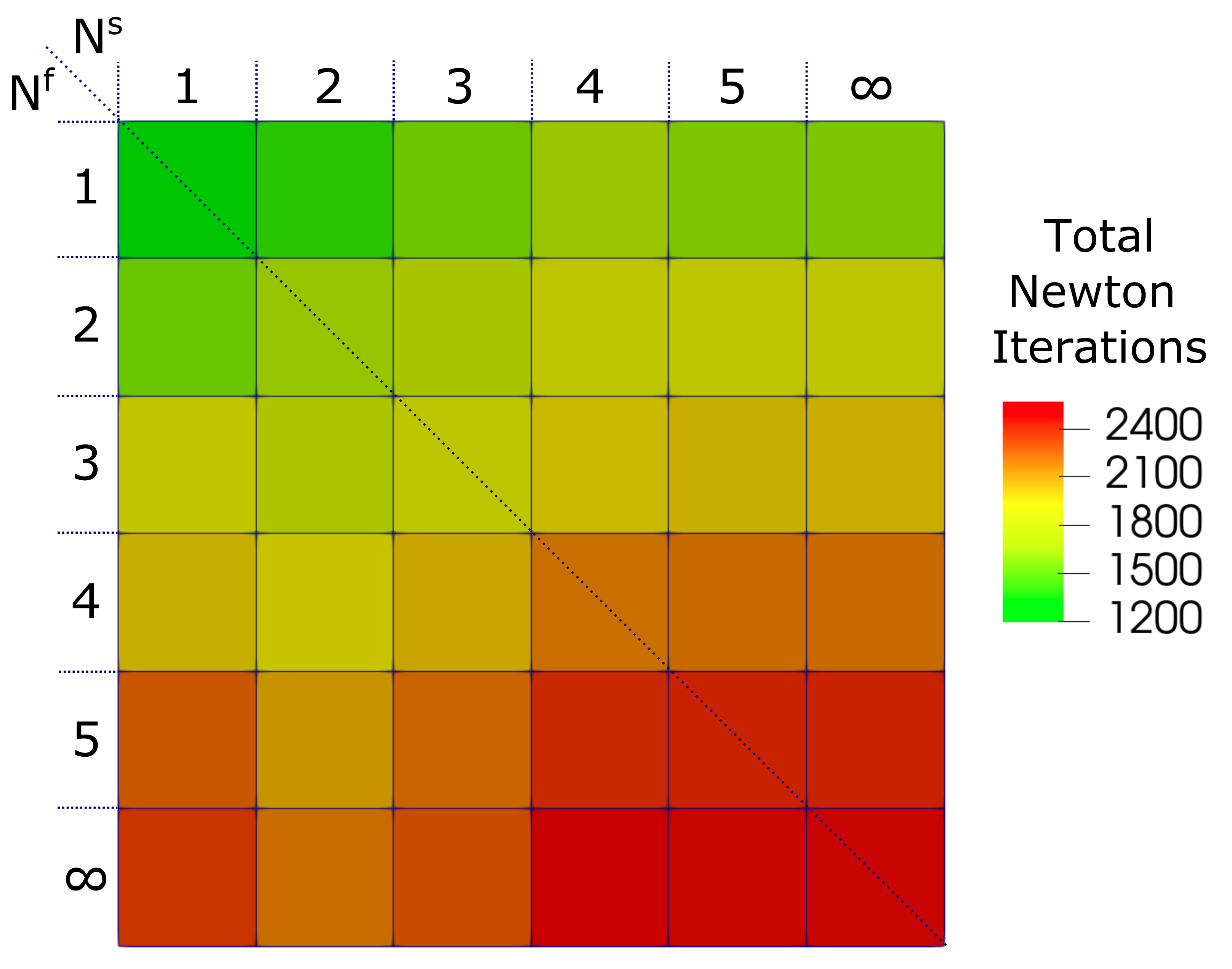}
	} 
	\caption{Influence of different choices for $N^f$ and $N^s$ on the required number of coupling $N_{Coupling}$ and Newton iterations $N_{Newton}$, illustrated by color scale.}
	\label{Fig:ContourPlots}
\end{figure}
Based on that, Table \ref{Tab:BottomCase-Newton} lists the number of coupling and Newton iterations -  for each subproblem as well as in total.
Figure \ref{Fig:ContourPlots} visualizes the data in two color plots.

As expected, the results show that increasing the number of Newton iterations per call
leads to fewer coupling iterations and vice versa.
Interestingly, however, the total number of Newton iterations $N_{Newton}$ shows the opposite trend: 
Running fewer or even just one Newton iteration per solver call 
requires fewer Newton steps.
Following the arguments from Section \ref{Sec:NewtonCost}, the computational cost are therefore expected to be lower too.

An explanation is that exchanging data  after every Newton iteration makes sure to always use the most recent solution of the other subproblem, keeping the boundary conditions up to date.
Moreover, in case an IQN method is applied, every Newton iteration adds a new data pair and hence improves 
the Jacobian approximation.

Another interesting observation is that
if the Newton iterations per solver call 
are kept fixed for one problem, increasing it for the other one reduces $N_{Coupling}$.
This effect can be very useful if the cost of the two subproblems are very different.
A common example would be a very complex fluid problem coupled to a rather coarse structural simulation:
In that case
setting $N^f=1$ and $N^s>1$ is expected to yield the best performance.

\subsection{Example: Elastic Beam in Channel Flow} \label{Sec:ChannelCase}

\begin{figure}[h!]
	\centering
	\subfloat[Test case geometry and material parameters.]
	{
		\resizebox{0.51\textwidth}{!}{	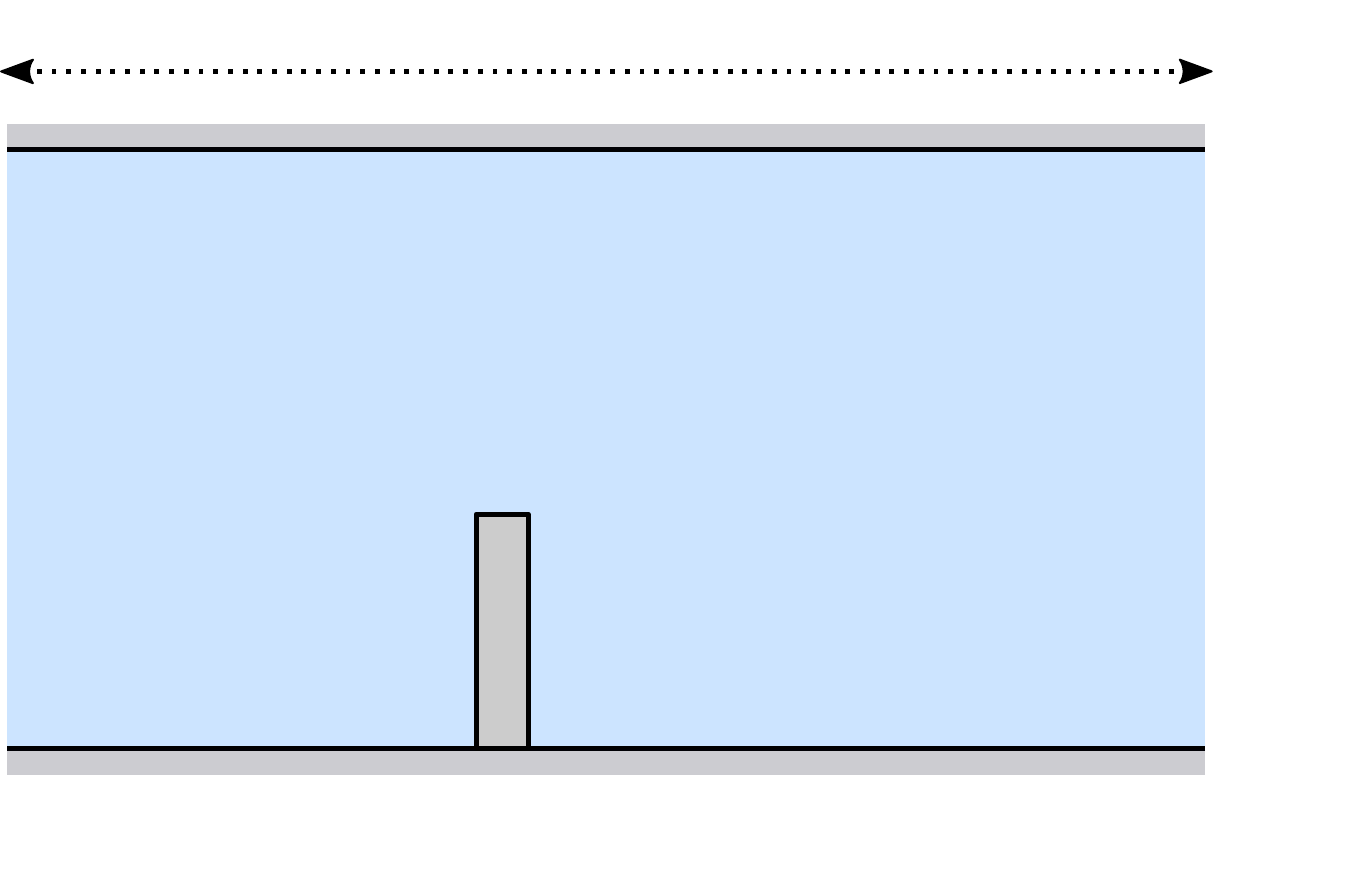}
	}  	
~~
	\subfloat[Snapshot in deformed state ($t=0.25\,\text{s}$).]
	{
		\includegraphics[trim=180 -20 300 0,clip, width=0.46\textwidth]{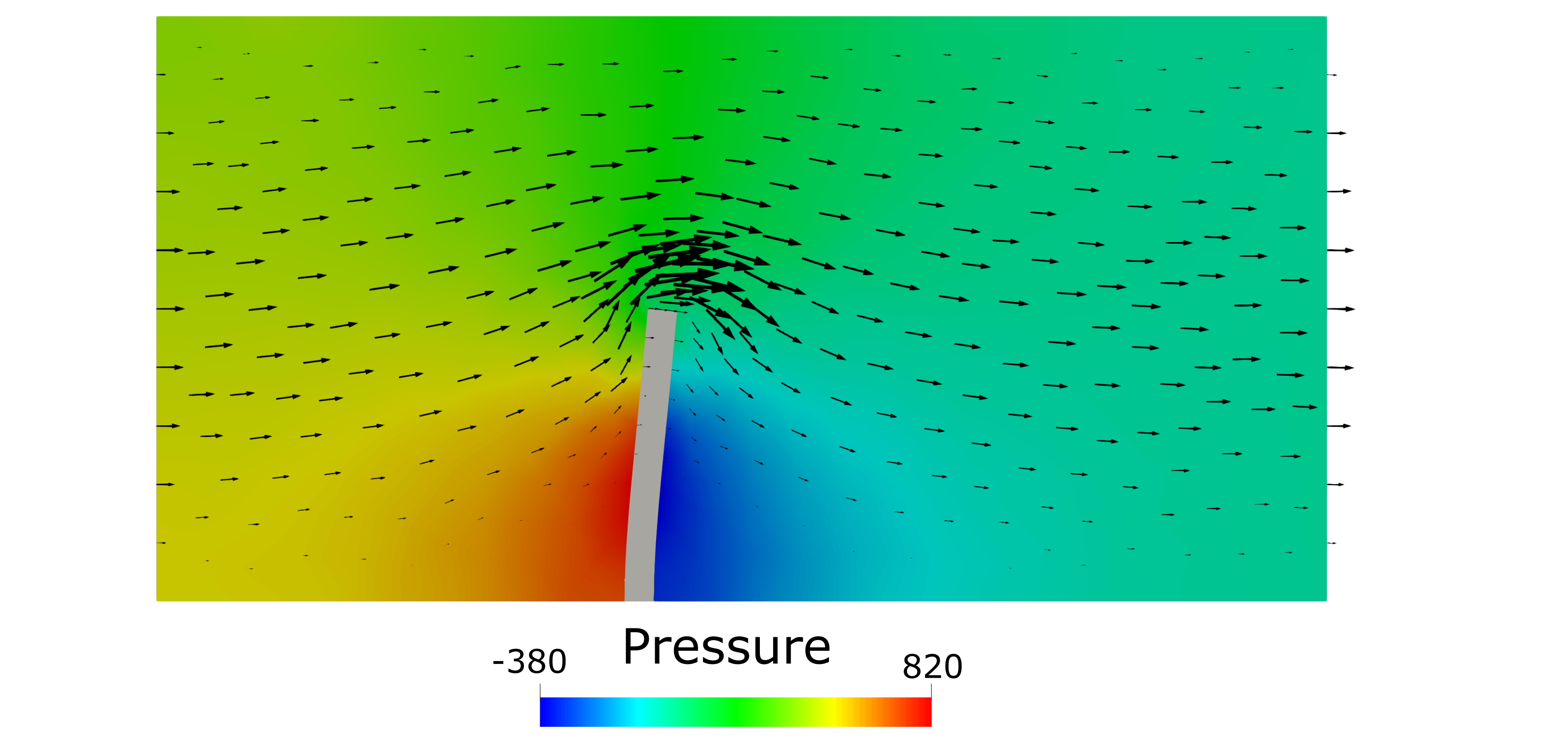}
	}  
	\caption{Elastic beam in channel flow test case.}
	\label{Fig:ChannelCase}
\end{figure}

While the observations made for the first example are typical for FSI problems with a strong interdependency,
the characteristics of the second test case, illustrated in Figure \ref{Fig:ChannelCase}, are different: Since the elastic beam positioned in the channel flow is rather heavy, the added-mass effect is less emphasized. Instead, the 
coupled system is mainly driven by the flow problem, so that an under-relaxation of the interface deformation with $\omega=0.8$ is sufficient to stabilize it.

The fluid problem is discretized by $771$ triangular finite elements, the structure by $60$ isogeometric elements of spline degree $2$.
The simulations run for $50$ time steps of width $\Delta t=0.005$.
Convergence is triggered by $\varepsilon_{Coupling}=10^{-5}$ and $\varepsilon_{Problem}=10^{-10}$. \\

\begin{figure}[h!]
	\centering
	\subfloat[Coupling iterations $N_{Coupling}$.]
	{
		\includegraphics[trim=0 0 0 0,clip, width=0.36\textwidth]{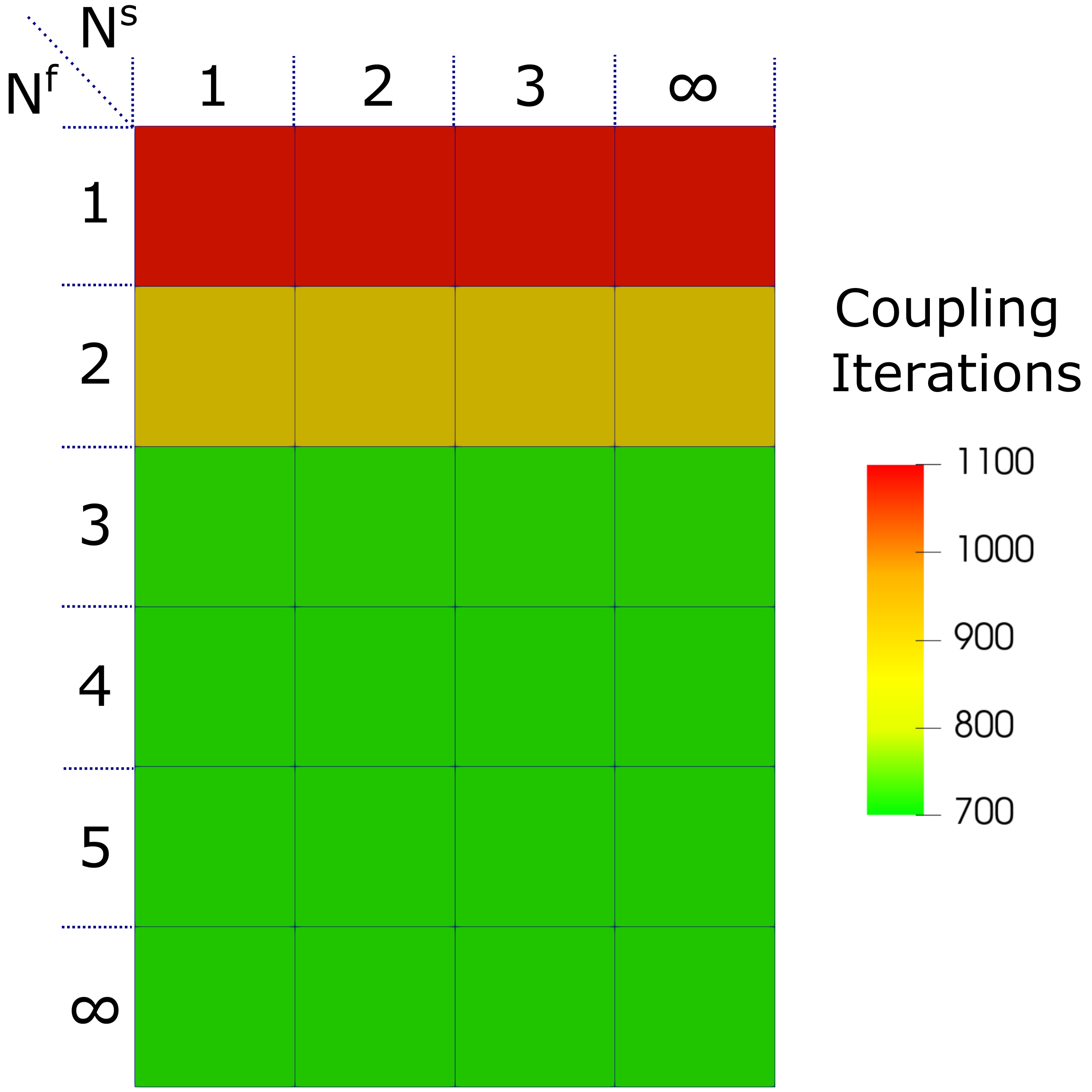}
	}  	
	\qquad
	\subfloat[Total Newton iterations $N_{Newton}$.]
	{
		\includegraphics[trim=-25 5 -25 0,clip, width=0.39\textwidth]{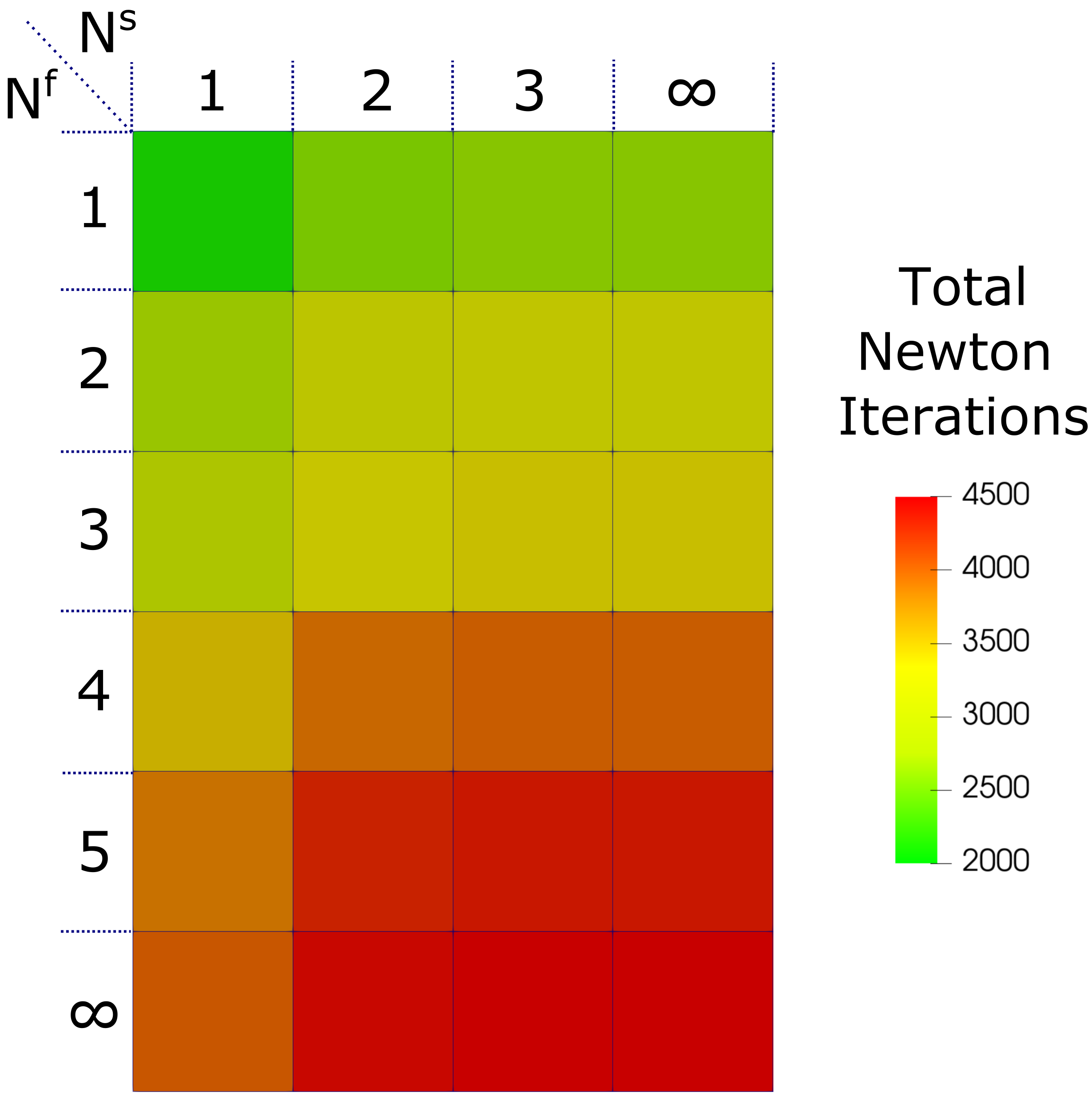}
	} 
	\caption{Influence of different choices for $N^f$ and $N^s$ on the required number of coupling $N_{Coupling}$ and Newton iterations $N_{Newton}$, illustrated by color scale.}
	\label{Fig:ContourPlotsCase2}
\end{figure}

Just as in the previous example, Table \ref{Tab:ChannelCase-Newton} lists the required coupling iterations $N_{Coupling}$, the single-field Newton iterations, as well as 
the total Newton iterations $N_{Newton}$, depending on the choices for $N^f$ and $N^s$, i.e., the Newton steps per solver call.
Figure \ref{Fig:ContourPlotsCase2} visualizes the data in two heat maps.
The columns for $N^s>3$ are omitted, as
none of the simulations required more than three Newton iterations 
in any structural solver call.

\begin{table}[h!]
		\caption{Newton and coupling iterations required for different choices of $N^f$ and $N^s$. The formatting is equivalent to that of Table \ref{Tab:BottomCase-Newton}: \textbf{coupling iterations}, \underline{total Newton iterations}, \textcolor{blue}{fluid's Newton iterations}, \textcolor{brown}{structural Newton iterations}.}
	\centering
	\small
	\begin{tabular}{ c   c  || c c | c  c | c c |  c c }
		%
		\multicolumn{2}{ c || }{ } &  \multicolumn{8}{ c  } { \textbf{Structural Newton Iterations $N^s$} } \\
		& & \multicolumn{2}{| c |}{\textbf{1}} & \multicolumn{2}{| c |}{\textbf{2}} & \multicolumn{2}{| c |}{\textbf{3}} &  \multicolumn{2}{| c }{$\mathbf{\infty}$} \\
				\hline \hline
		\multirow{1}{*}{ \rotatebox{90}{ \textbf{Fluid Newton Iterations $N^f$} }}
		& \multirow{2}{*}{\textbf{1}}    &  \iterC{1083} & \iterN{2109} & \iterC{1083} & \iterN{2548}  &  \iterC{1083} & \iterN{2611}  &  \iterC{1083} & \iterN{2611}  \\
		& 	&   \iterF{1083} &  \iterS{$1026$} & \iterF{1083} &  \iterS{$1465$} & \iterF{1083} &  \iterS{$1528$} & \iterF{1083} &  \iterS{$1528$}  \\
		\cline{2-10}
		& \multirow{2}{*}{\textbf{2}}    &  \iterC{901} & \iterN{2702} & \iterC{901} & \iterN{3142}  &  \iterC{901} & \iterN{3205}  &  \iterC{901} & \iterN{3205}  \\
		& 	&   \iterF{$1802$} &  \iterS{$900$} & \iterF{$1802$} &  \iterS{$1340$} & \iterF{1802} &  \iterS{$1403$} & \iterF{1802} &  \iterS{$1403$} \\
		\cline{2-10}
		& \multirow{2}{*}{\textbf{3}}    &  \iterC{721} & \iterN{2884} & \iterC{721} & \iterN{3324}  &  \iterC{721} & \iterN{3386}  &  \iterC{721} & \iterN{3386}  \\
		& 	&   \iterF{2163} &  \iterS{721} & \iterF{2163} &  \iterS{1161} & \iterF{2163} &  \iterS{1223} & \iterF{2163} &  \iterS{1223} \\
		\cline{2-10}
		& \multirow{2}{*}{\textbf{4}}    &  \iterC{718} & \iterN{3485} & \iterC{718} & \iterN{3925}  &  \iterC{718} & \iterN{3988}  &  \iterC{718} & \iterN{3988}  \\
		& 	&   \iterF{2767} &  \iterS{718} & \iterF{2767} &  \iterS{1158} & \iterF{2767} &  \iterS{1221} & \iterF{2767} &  \iterS{1221}  \\
		\cline{2-10}
		& \multirow{2}{*}{\textbf{5}}    &  \iterC{718} & \iterN{3867} & \iterC{718} & \iterN{4307}  &  \iterC{718} & \iterN{4370}  &  \iterC{718} & \iterN{4370}  \\
		& 	&   \iterF{3149} &  \iterS{718} & \iterF{3149} &  \iterS{1158} & \iterF{3149} &  \iterS{1221} & \iterF{3149} &  \iterS{1221}  \\
		\cline{2-10}
		%
		& \multirow{2}{*}{$\mathbf{\infty}$}    &  \iterC{718} & \iterN{4014} & \iterC{718} & \iterN{4454}  &  \iterC{718} & \iterN{4517}  &  \iterC{718} & \iterN{4517}  \\
		& 	&   \iterF{3296} &  \iterS{718} & \iterF{3296} &  \iterS{1158} & \iterF{3296} &  \iterS{1221} & \iterF{3296} &  \iterS{1221} \\
	\end{tabular}
	\label{Tab:ChannelCase-Newton}
\end{table}

This already points out the biggest difference to Table \ref{Tab:BottomCase-Newton},
the negligible influence of $N^s$ on the coupling iterations $N_{Coupling}$.
It is a result of the structure's lower sensitivity to the flow solution, i.a., due to the decreased density ratio of $\rho^f / \rho^s = 0.1$.

For the number of Newton iterations performed per fluid solver call $N^f$, in contrast, the same relation to $N_{Coupling}$ as in the previous test case can be observed:
While increasing $N^f$ results in fewer coupling steps, 
choosing $N^f=N^s=1$ again yields the smallest number of Newton iterations and therefore supposedly also the 
lowest computational cost.

\subsection{Adapt Newton Iterations Dynamically}  \label{Sec:AdaptiveMethods}

The results discussed so far identified running just one Newton iteration per solver call as a good and efficient choice.
Nevertheless,
it is not without flaws in that it typically requires many coupling iterations, which
might produce unwanted coupling overhead.
This motivated the idea of adapting $N^f$ and $N^s$ dynamically.
We propose two approaches:
\begin{enumerate}
	\item The frequent data exchange is of particular importance in the first Newton iterations, in which the solutions still change a lot. In later steps, however, the increments are decreasing until they might no longer justify the effort of coupling after every iteration.
	As a remedy, we suggest to run just one Newton step per solver call, $N^f=N^s=1$, as long as coupling convergence (see Section \ref{Sec:Convergence}) is not reached. Once it is, the two solvers switch to full single-field convergence, $N^f=N^s=\infty$. 
	We will refer to this technique as \textit{N1-CC} \textit{approach}  (``one Newton iteration until coupling convergence").
	Naturally, further variants are obtained by running more than one iteration for the unconverged coupling, i.e., N2-CC, N3-CC, etc.\\
	Note that in case after setting $N^f=N^s=\infty$ the coupling convergence criterion is not satisfied anymore,
	$N^f$ and $N^s$ are switched back. In practice, however, this issue never occurred in any simulation of this work.
	If it would, using a stricter criterion, e.g., $0.1 \cdot \varepsilon_{Coupling}$, for triggering full single-field convergence could fix it.
  	\item The goal of the second suggestion is to avoid feeding back inaccurate data into the coupling loop. Therefore,
  	the \textit{converged interface data} approach runs the Newton loop until the relative change in the coupling data, e.g., the fluid tractions, is less than some bound $\varepsilon_{CID}$. 
  	In contrast to the N1-CC method, this approach typically results in additional Newton steps only in the first few coupling iterations.
\end{enumerate}

Revisiting the two test cases
from Sections \ref{Sec:BottomCase} and \ref{Sec:ChannelCase}, 
Table \ref{Tab:AdaptiveMethods} investigates the effectiveness of these adaptive techniqes.
\begin{table}[h!]
	\caption{Iterations required for the adaptive techniques (formatting as in previous tables).
		The results in gray are repeated for simpler comparison.}
\centering
	\small
	\begin{tabular}{  c  || c c | c  c }
		%
		\textbf{Setting} &  \multicolumn{2}{ c  } { \textbf{Elastic Bottom}} & \multicolumn{2}{ | c } { \textbf{Beam in Channel Flow}} \\
				\hline \hline
	\multirow{2}{*}{ \textbf{N1-CC}} & \iterC{417} & \iterN{945} & \iterC{768} & \iterN{1637}  \\
	& \iterF{528}  & \iterS{417}  & \iterF{869}  &\iterS{768}  \\
	\hline 
		\multirow{2}{*}{ \textbf{N3-CC}} & \iterC{355} & \iterN{1799} & \iterC{720} & \iterN{3384}  \\
	& \iterF{1071}  & \iterS{728}  & \iterF{2161}  &\iterS{1223}  \\
	\hline 
	\multirow{2}{*}{ \textbf{Converged Interface Data}, $\varepsilon_{CID}=10^{-4}$ } &  \iterC{572} & \iterN{1720} & \iterC{1083} & \iterN{3182}  \\
	& \iterF{928}  & \iterS{792}  & \iterF{1880}  &\iterS{1302}  \\
\hline \hline
	\multirow{2}{*}{ \textcolor{gray}{$N^f=N^s=1$}} &  \iterCompare{617} &  \iterCompare{1166} &  \iterCompare{1083} &  \iterCompare{2109}  \\
&  \iterCompare{617}  &  \iterCompare{567}  & \iterCompare{1083}  &  \iterCompare{1026}  \\
\hline 
\multirow{2}{*}{ \textcolor{gray}{$N^f=N^s=\infty$}} &  \iterCompare{354} &  \iterCompare{2547} &  \iterCompare{718} &  \iterCompare{4517}  \\
&  \iterCompare{1700}  &  \iterCompare{847}  & \iterCompare{3296} &  \iterCompare{1221}  \\
	\end{tabular}
\label{Tab:AdaptiveMethods}
\end{table}
It clearly shows that the suggested N1-CC variant significantly reduces the number of coupling iterations $N_{Coupling}$ compared to the choice 
$N^f=N^s=1$, 
overcoming its major drawback. 
At the same time, its main strength, i.e., the very low
total number of Newton iterations $N_{Newton}$, is retained.
In fact, the results indicate that $N_{Newton}$ is even further decreased for the given examples. 

These observations can be explained by the reduced data exchange of the N1-CC approach
once coupling convergence is reached:
Aside from the evident effect on $N_{Coupling}$,
the Newton iterations are no longer synchronized.
Therefore, an already converged problem will not be called for each iteration
of the other, non-converged problem, avoiding unnecessary overhead.
Moreover, depending on $\varepsilon_{Coupling}$ and $\varepsilon_{Problem}$,
the remaining increments might satisfy the coupling convergence,
yet still
slow down the other problem's Newton loop by slightly updating its boundary conditions in every iteration.

Of course these investigations are mainly based on examples, but nevertheless they indicate a superiority of the dynamic N1-CC method over
using fixed numbers of iterations.
Beyond that, the results of the N3-CC variant are interesting for cases in which more focus is put on the coupling iterations, e.g., if the data exchange is expensive,
because it reaches the same number of coupling iterations as $N^f=N^s=\infty$
in less Newton steps.
The converged interface data approach, 
on the other hand, cannot keep up with neither $N^f=N^s=1$ nor the N1-CC approach. 
This indicates that the risk of passing back inaccurate data into the coupling loop is not very severe.
Instead, the additional Newton steps in the first coupling iterations, in which the interface data has
not converged yet, 
prove to be rather ineffective,
which is in line with the discussion on the N1-CC approach.

\section{Conclusion}

This work revolves around one central question:
In partitioned fluid-structure interaction with two non-linear subproblems, 
how many Newton iterations per solver call result in the most efficient coupling scheme?

Its motivation arises from the argument 
that the required Newton iterations are a better measure for a partitioned algorithm's computational efficiency
than the number of coupling iterations;
because rather than assuming constant cost per solver call,
they reflect the reiterative nature of the numerical solution procedure, i.e., the Newton loop.
Based on typical examples, the discussion shows that iterating to full convergence for every solver call does in fact
require the fewest coupling steps. The total number of Newton iterations and therefore 
the computational cost, however, can be reduced significantly
by a more frequent communication, in particular by running just one Newton step per call.

Against this backdrop, this work discusses two adaptive choices for the number of Newton iterations per solver call.
In particular, the N1-CC approach exchanges data after every Newton iteration only until coupling convergence is reached; 
after that, every solver call iterates to full convergence.
The numerical examples confirm that this technique maintains (and even further decreases) the low number of Newton iterations,
while 
reducing the number of coupling iterations substantially.

Although the discussion is purely based on numerical experiments and logical arguments rather than a firm mathematical foundation, 
its
findings provide valuable assistance and guidelines 
on how to properly set the
number of 
Newton iterations per solver call 
in a partitioned fluid-structure interaction scheme.

In the end, the optimal choice of course depends on a variety of factors, as for example the time step size, the specific implementation, or in general the problem at hand.
Therefore, future works will extend the investigation to a wider set of parameters.


\begin{thebibliography}{99}

\bibitem{lindner2015comparison} Lindner, F. and Mehl, M. and Scheufele, K. and Uekermann, B. A comparison of various quasi-Newton schemes for partitioned fluid-structure interaction. \textit{COUPLED VI: proceedings of the VI International Conference on Computational Methods for Coupled Problems in Science and Engineering} (2015) 477--485.

\bibitem{Degroote2009} Degroote, J. and Bathe , K.J. and Vierendeels, J. Performance of a new partitioned procedure versus a monolithic procedure in fluid--structure interaction. \textit{Computers {\&} Structures} (2009) \textbf{87}:793--801.

\bibitem{spenke2020multi} Spenke, T. and Hosters, N. and Behr, M. A multi-vector interface quasi-Newton method with linear complexity for partitioned fluid--structure interaction. \textit{Comput.  Methods Appl.  Mech.  Eng.} (2020) \textbf{361}:112810.

\bibitem{pauli2016stabilized} Pauli, L. and Behr, M. On stabilized space-time FEM for anisotropic meshes: Incompressible Navier--Stokes equations and applications to blood flow in medical devices. 
\textit{Int. J. Numer. Methods Fluids} (2017) \textbf{85}:189--209.

\bibitem{donea2003finite} Donea, J. and Huerta, A. \textit{Finite element methods for flow problems}. John Wiley \& Sons, (2003).

\bibitem{forti2015semi} Forti, D. and Ded{\`e}, L. Semi-implicit BDF time discretization of the Navier--Stokes equations with VMS-LES modeling in a high performance computing framework. \textit{Computers \& Fluids} (2015) \textbf{117}:168--182.

\bibitem{emumPaper} Behr, M. and Abraham, F. Free-surface flow simulation in the presence of inclined walls, \textit{Comput.  Methods Appl.  Mech.  Eng.} (2002) \textbf{191}:5467--5483.

\bibitem{bathe2006finite} Bathe, K.J. \textit{Finite element procedures}. TBS, (1996).

\bibitem{yibin2001nonlinear} Yibin, F.U. and Ogden, R.W. \textit{Nonlinear elasticity: Theory and applications}, London Mathematical Society Lecture Note Series, Vol. 283, (2001).

\bibitem{cottrell2009isogeometric} Cottrell, J.A. and Hughes, T.J.R. and Bazilevs, Y. \textit{Isogeometric analysis: Toward integration of CAD and FEA}. John Wiley \& Sons, (2009).

\bibitem{Hughes} Hughes, T.J.R. and Cottrell, J.A. and Bazilevs, Y. Isogeometric analysis: CAD, finite elements, NURBS, exact geometry and mesh refinement. \textit{Comput.  Methods Appl.  Mech.  Eng.} (2005) \textbf{194}:4135--4195.

\bibitem{chung1993time} Chung, J. and Hulbert, G.M. A time integration algorithm for structural dynamics
with improved numerical dissipation: The generalized-$\alpha$ method. \textit{Journal of Applied Mechanics} (1993) \textbf{60}:371–375.

\bibitem{erlicher2002analysis} Erlicher, S. and Bonaventura, L. and Bursi, O.S. The analysis of the generalized-$\alpha$ method for non-linear dynamic problems. \textit{Computational mechanics} (2002) \textbf{28}:83--104.

\bibitem{Hosters2017} Hosters, N. and Helmig, J. and Stavrev, A. and Behr, M. and Elgeti, S. Fluid-structure interaction with NURBS-based coupling. \textit{Comput.  Methods Appl.  Mech.  Eng.} (2018) \textbf{332}:520--539.

\bibitem{kuttler2006solution} K{\"u}ttler, U. and F{\"o}rster, C. and Wall, W.A. A solution for the incompressibility dilemma in partitioned fluid--structure interaction with pure Dirichlet fluid domains. \textit{Computational Mechanics} (2006) \textbf{38}:417--429.

\bibitem{degroote2011multi} Degroote, J. and Vierendeels, J. Multi-solver algorithms for the partitioned simulation of fluid--structure interaction. \textit{Comput.  Methods Appl.  Mech.  Eng.} (2011) \textbf{200}:2195--2210.

\bibitem{hostersspline} Hosters, N. \textit{Spline-based methods for fluid-structure interaction}. Dissertation, RWTH Aachen University, (2018).

\bibitem{forster2007robust} F{\"o}rster, C. \textit{Robust methods for fluid-structure interaction with stabilised finite elements}. Dissertation, University of Stuttgart, (2007).

\bibitem{forster2007artificial} F{\"o}rster, C. and Wall, W.A. and Ramm, E. Artificial added mass instabilities in sequential staggered coupling of nonlinear structures and incompressible viscous flows. \textit{Comput.  Methods Appl.  Mech.  Eng.} (2007) \textbf{196}:1278--1293.

\bibitem{Causin} Causin, P. and Gerbeau, J.-F. and Nobile, F. Added-mass effect in the design of partitioned algorithms for fluid-structure problems. \textit{Comput.  Methods Appl.  Mech.  Eng.} (2005) \textbf{194}:4506--4527.

\bibitem{kuttler2008fixed} K{\"u}ttler, U. and Wall, W.A. Fixed-point fluid--structure interaction solvers with dynamic relaxation. \textit{Computational mechanics} (2008) \textbf{43}:61--72.

\bibitem{irons1969version} Irons, B.M. and Tuck, R.C. A version of the Aitken accelerator for computer iteration. \textit{Int. J. Num. Meth. Engng.} (1969) \textbf{1}:275--277.

\end{thebibliography}
\end{document}